\documentclass{amsart}
\usepackage{amssymb}
\usepackage{amsfonts}
\usepackage{amsmath}
\usepackage{lastpage}
\usepackage[legalpaper,bookmarks=true,colorlinks=true,linkcolor=blue,citecolor=blue]{hyperref}
\usepackage{graphicx}%
\usepackage{fancyhdr}
\usepackage{color}
\usepackage[mathlines]{lineno}
\usepackage{lscape}
\usepackage{epsfig}
\usepackage{natbib}
\usepackage{geometry}
\usepackage{tgbonum}
\fontfamily{qcr}\selectfont
\usepackage[]{listings}

\newtheorem{theorem}{Theorem}
\theoremstyle{plain}

\newtheorem{corollary}{Corollary}

\numberwithin{equation}{section}

\newcommand{\Bin}{\bigskip \noindent}
\newcommand{\Bi}{\bigskip}
\newcommand{\Ni}{\noindent}

\begin{document}
\Large
\title[]{Extremes, extremal index estimation, records, moment problem for the Pseudo-Lindley distribution and applications}

\begin{abstract} The pseudo-Lindley distribution which was introduced in Zeghdoudi and Nedjar (2016) is studied with regards to its upper tail. In that  regard, and  when the underlying distribution function follows the Pseudo-Lindley law, we investigate  the  behavior of its values, the asymptotic normality of the Hill estimator and the double-indexed generalized Hill statistic process (Ngom and Lo), the asymptotic normality of the records values and the moment problem.\\

\noindent $^{\dag}$ Gane Samb LO.\\
LERSTAD, Gaston Berger University, Saint-Louis, S\'en\'egal (main affiliation).\newline
LSTA, Pierre and Marie Curie University, Paris VI, France.\newline
AUST - African University of Sciences and Technology, Abuja, Nigeria\\
gane-samb.lo@edu.ugb.sn, gslo@aust.edu.ng, ganesamblo@ganesamblo.net\\
Permanent address : 1178 Evanston Dr NW T3P 0J9,Calgary, Alberta, Canada.\\

\noindent \noindent $^{\dag}$ Dr. Modou NGOM\\
Work Affiliation : Ministery of High School (SENEGAL)\\
LERSTAD, Gaston Berger University, Saint-Louis, S\'en\'egal.\newline
ngomodoungom@gmail.com\\

\noindent $^{\dag}$ Dr. Moumouni Diallo\\
Université des Sciences Sociale et de Gestion de Bamako ( USSGB)\\
Faculté des Sciences Économiques et de Gestion (FSEG)\\
Email: moudiallo1@gmail.com.\\

\end{abstract}

\maketitle

\section{Introduction}\label{sec1}

\Ni \textbf{1. General facts}.\\

\noindent The following probability distribution function (\textit{pdf}), named as the Pseudo-Lindley \textit{pdf}, 

\begin{equation} 
f(x)=f(x,\theta,\beta)=\frac{\theta(\beta-1+\theta x) e^{-\theta x}}{\beta} 1_{(x\geq 0)} \label{glind}
\end{equation}

\Bin with parameters $\theta>0$ and $\beta>1$, has been introduced by \cite{Zeghoudi3}  as a generalization of the Lindley \textit{pdf} :

\begin{equation} 
\ell(x)=\frac{\theta^2(1+x) e^{-\theta x}}{1+\theta} 1_{(x\geq 0)}. \label{olind}
\end{equation}

\Bin  in the sense that for $\beta=1+\theta$, $f(\circ)$ is identical to $\ell(\circ)$.\\

\Ni  Actually, $f$ derives from $\ell$ by a mixture of a Lindley distributed random variable and an independent $\Gamma(2,\theta)$ random variables with mixture coefficients $r_1=(\beta-1)/\beta$ and $r_2=1/\beta$, where $1<r_1, \ r_2<1$ and $r_1+r_2=1$.\\

\Ni The cumulative distribution \textit{cdf} function is given by

$$
1-F(x)=\left(\beta^{-1}(\beta+ \theta x) e^{-\theta x}\right) 1_{(x\geq 0)}.
$$

\Bin The Lindley original distribution is an important law that been used and still being used in Reliability, in Survival analysis and other important disciplines. Because of its original remarkable qualities, it kicked off a considerable number generalizations as pointed out by \cite{Zeghoudi3}. The current generalization
\eqref{glind} has been tested on real data and simulated and shows real interest in survival analysis, on the Guinean Ebola data for example (\cite{Zeghoudi3}). The paper of focused on asymptotic tests of that law based on moments estimators. The interest that distribution demonstrated in real data modeling motivated us to give some asymptotic theories on it, in view of statistical tests. In this paper, we deal with the properties of the upper tail, the extreme value distribution and the record values. etc., which of them providing statistical tests.\\

\noindent Throughout the paper, $X$, $X_1$, $X_2$, $\cdots$ is a sequence independent real-valued random variables (\textit{rv}), defined on the same probability space $(\Omega, \mathcal{A},\mathbb{P})$, with common cumulative distribution function $F$, with the first asymptotic moment function and the generalized inverse function  defined by

$$
R(x,F)=\frac{1}{1-F(x)}\int_{x}^{+\infty}(1-F(y)) \ dy, \ x \in ]0,+\infty[
$$

\bigskip \noindent and

$$
F^{-1}(u)=\inf \{x \in \mathbb{R}, \ F(x)\geq u\} \ for \ u \in ]0,1[ \ and \ F^{-1}(0)=F^{-1}(0+).
$$

\Bin For each $n\geq 1$, we denote the ordered statistics of the sample $X_1$, $\cdots$, $X_n$ by 

$$
X_{1,n}\leq \cdots \leq X_{n,n}.
$$ 

\Bin Usually, in extreme value theory, we focus on upper extreme and the hypothesis $X>0$ and the log-transform $Y=\log X$ is instrumental in all major results in that field. We denote $G(x)=F(e^x)$, $x \in \mathbb{R}_{+}$. The Renyi representation is also of common use in the following form. The sequence is replaced as follows

\begin{equation}
\{ \{X_{1,n}\leq \cdots \leq X_{n,n}\}, \ n\geq 1\}=_{d} \{ \{F^{-1}(1-U_{n-j+1,n}), 1\leq j\leq n\}, \ n\geq 1\}, \label{renyi}
\end{equation}

\Bin where  $=_{d}$ stands for the equality in distribution. Finally, the following Malmquist representation (see \cite{shwell}, also \cite{ips-wcia-ang}, page ...) is also used : for each $n\geq 1$, there exist a finite sequence of standard independent exponential random variables $E_{1,n}$, $\cdots$, $E_{n,n}$ such that

\begin{equation} \label{malmquist}
\left\{ \left(\frac{U_{i+1,n}}{U_{i,n}}\right)^{i}, \ 1\leq i \leq n \right\}=_{d} \left\{E_{i,n}, 1\leq i \leq n \right\}.
\end{equation}

\Bi

\section{Extremes}

\noindent We can directly see that $F$ is the Gumbel distribution $G_{0}$ by three different arguments. First, by using the Von Mises' argument (see \cite{dehaan} or
\cite{ips-wciia-ang}, Proposition 24, page 184)

\begin{eqnarray}
\lim_{x \rightarrow 0} \frac{f^{\prime}(x)(1-F(x))}{f^2(x)}=-1. 
\end{eqnarray}

\Bin A second argument comes from that $Y=\exp(X)$ has the distribution $G(x)=F(\log x)=\beta^{-1} (\beta + \theta \log x) x^{-\theta x}$. Since

\begin{eqnarray}
\lim_{x \rightarrow 0} \frac{1-G(\lambda x)}{1-G(x)}=\lambda^{-\theta}, 
\end{eqnarray}

\Bin $G \in D(G_{1/\theta})$ and since $F(x)=G(e^x)$ for $x\geq 1$,  by Theorem (in Lo), $F \in D(G_0)$.\\

\Ni A third argument is relation to the development of the quantile function. In the appendix (page \pageref{appendix-erZ}), we give a number of expansions of that quantile that could be used for different purposes. For example we have (see page \pageref{erZ_12d}),

\begin{equation} \label{quantile}
F^{-1}(1-u) = \theta^{-1}(\log(1/u) -  \log\log(1/u)) +  \theta^{-1} K(u)
\end{equation}
 
\Bin with $K(u)=O(\log 1/u)^{-2})$. By using it we get 

$$
\frac{F^{-1}(1-\lambda u)-F^{-1}(1-u)}{(1/\theta)} \rightarrow -\log \lambda \ as \ u\rightarrow 0.
$$

\Ni By the $\pi$-variation criteria of \cite{dehaan} (See \cite{ips-wcia-ang}, Proposition 11, page 88), we have $F\in D(G_0)$ and $R(x,F)\rightarrow \gamma=1/\theta$ as $x\rightarrow +\infty$. Formula \eqref{quantile} is actually a second-order condition for the quantile function (see  \cite{dehaan}). We apply it right to get a rate of convergence of the maximum observations. Put $\gamma=1/\theta$.

\Bi \textbf{2. Expansion of the maximum values}.\\

\Ni By the Renyi representation and by denoting  $Z_n=-\log(nU_{1,n})$, we have that $\log(1+Z_n/(\log n)) \rightarrow_{\mathbb{P}} 0$ and since $\log U_{1,n}=O_{\mathbb{P}}(\log n)^{-1}$

\begin{eqnarray*} 
X_{n,n}-F^{-1}(1-1/n)&=&\gamma Z_n + \gamma \log(1+Z_n/(\log n))\\
&+& O((\log n)^{-2})+O((\log U_{1,n})^{-2}).
\end{eqnarray*}

\noindent and hence

\begin{eqnarray} 
\frac{X_{n,n}-F^{-1}(1-1/n)}{\gamma}= Z_n + O_{\mathbb{P}}(\log n)^{-1})=\Gamma + O_{\mathbb{P}}(1)  
\end{eqnarray}

\Bin It is easy to see that $Z_n$ converge to Gumbel law $\Lambda$ with \textit{cdf} 

$$
G_{0}(x)=\exp(-\exp(-x)), \ x \in \mathbb{R}
$$

\Bin So we have that $X_{n,n}$ converges to a $\Lambda$ law. But we obtain the random rate of convergence $Z_n/\log n$ since

$$
\frac{\log Z_n}{\log n} \left(\frac{X_{n,n}-F^{-1}(1-1/n)}{\gamma}-Z_n \right)=1
$$

\Bin As well for $k=k(n)\rightarrow +\infty$ such that $k(n)/n\rightarrow 0$, and by taking $T_n=\log (nU_{k,n}/k)$ and $q_n=n/k(n)$ which goes to $+\infty$, we have

\begin{eqnarray} 
\frac{X_{n-k,n}-F^{-1}(1-k/n)}{\gamma}= T_n + \log(1+T_n/\log q_n)) + O_{\mathbb{P}}((\log q_n)^{-2}).
\end{eqnarray}

\Bin \textbf{3. Estimating the extreme value index  $\gamma=1/\theta$}.\\

\noindent The \cite{hill}'s estimator 

\begin{equation}
H_n=\frac{1}{k(n)} \sum_{j=1}^{k(n)} j\left(X_{n-j+1,n}-X_{n-j,n}\right),
\end{equation}

\Bin is the most celebrated estimator the extreme value index $\gamma=1/\theta$ of $Z=\exp(X)$. Among a significant number of generalizations, the \cite{ngomlo2016}'s generalization, called the Double-indexed function Hill estimator, is one the sharpest one. It is defined as

$$
H_n(f,s)= \left(\sum_{j=1}^{k(n)} f(j)\left(X_{n-j+1,n}-X_{n-j,n}\right)^s/(a_n(f,s)\right)^{1/s},
$$

\noindent where $f:\mathbb{N}\setminus \{0\} \rightarrow \mathbb{R}_+\setminus \{0\}$ is a measurable mapping and $s>0$. Let us define

$$
a_n(f,s)=\Gamma(s+1) \sum_{j=1}^{k(n)} f(j)j^{-s}, \ C^2(s)=\Gamma(2s+1)-\Gamma(s+1)^2, \ s_n^2(f,s)= C^2(s) \sum_{j=1}^{k(n)} f(j)^2j^{-2s},
$$

\Bin and 
 
$$
B_n(f,s)=\max\{ f(j)j^{-s}/s_n(f,s), \ 1\leq j \leq k(n)\}.
$$

\Bin We simply notice that the classical Hill's estimator is $H_n(I_d,1)$ where $I_d$ is the identity function on $\mathbb{N}\setminus \{0\}$. Let us give  asymptotic normality for Double-indexed function Hill estimator.\\

\noindent \textbf{(a) Extreme Limit Theorem}.\\

\noindent We begin with the simple Hill's estimator.\\ 

\begin{theorem} \label{theoHill} For $]0, \ n] \ni k(n)\rightarrow +\infty$ such that 

\begin{equation}
k(n)^{3/4}/\log n \rightarrow 0. \tag{K1}
\end{equation}

\Bin we have, as $n\rightarrow +\infty$,

\begin{equation}
\sqrt{k(n)} \left(H_n-\gamma\right) \rightsquigarrow \mathcal{N}(0,\gamma^2). \label{hillSN}
\end{equation}

\end{theorem}

\Bin We want to establish the random rate of convergence associated with the convergence \ref{hillSN} in the part (a) of the following corollary. In the part (b), we want to share that we need any other condition on top of $k(n)/n \rightarrow 0$ to have the central limit theorem if $F^{-1}$ is reduced to 

\begin{equation}
F_{\ast}^{-1}(1-u)=\gamma \log u - C(\gamma) \log \log(1/u), \ u \in ]0,1[, \ C(\gamma)\geq 1. \label{Freduced}
\end{equation}   

\begin{corollary} \label{theoHillCoro} We have the following results.\\

\noindent (a) Here again $F$ is the \textit{cdf} of the Pseudo-Lindley distribution with parameters $\theta>0$ and $\beta>0$ and the notation above. Let $k(n)/\log n\rightarrow 0$.

\Bin Let $W(1)$ is a standard Gaussian random variable. Then we have 

$$
 \frac{\log n}{\gamma \sqrt{k(n)}} \biggr(\sqrt{k(n)}(H_n-\gamma)-\gamma W(1)\biggr) \rightarrow_{\mathbb{P}} 1,
$$

\Bin (b) If $F^{-1}$ were reduced as in Formula \eqref{Freduced}, we have the asymptotic normality 

$$
\sqrt{k(n)}(H_n-\gamma) \rightarrow \mathcal{N}(0,\gamma^2)
$$

\Bin whenever $k(n)/n \rightarrow 0$ and

$$
\log n \biggr(\sqrt{k(n)}(H_n-\gamma) - \gamma W(1)\biggr)=O_{\mathbb{P}}(1). \ \Diamond
$$
\end{corollary}

\Bin \textbf{Proof of Theorem \ref{theoHill}}. By Malmquist representations (See \cite{shwell} or \cite{ips-wcia-ang}, Proposition 32, page 135), by Formula \eqref{erZ_12e}, we have for any $1\leq j \leq k$,

\begin{eqnarray} 
X_{n-j+1,n}-X_{n-j,n}&=&F^{-1}(1-U_{j,n})-F^{-1}(1-U_{j+1,n}) \notag\\
&=& \gamma j^{-1} E_{j,n} - \gamma \int_{U_{j,n}}^{U_{j+1,n}} \frac{du}{u \log (1/u)} + O_{\mathbb{P}}\left(\left(log n\right)^{-2}\right) \label{spacingF}
\end{eqnarray}

\Bin and next 

\begin{eqnarray*}
j \left(X_{n-j+1,n}-X_{n-j,n}\right)&=& \gamma E_{j,n} - \gamma j\int_{U_{j,n}}^{U_{j+1,n}} \frac{du}{u \log (1/u)}  + O_{\mathbb{P}}\left(k \left(log n\right)^{-2}\right).
\end{eqnarray*}

\Bin So for $Z_n=\log nU_{1,n}$ (which converges in law to  $\Lambda$) and

\begin{equation} \label{step_01}
\left| \int_{U_{j,n}}^{U_{j+1,n}} \frac{du}{u \log (1/u)} \right|\leq \frac{j^{-1} E_{j,n}}{|\log n - Z_n|}.
\end{equation}

\Bin Hence

$$
\frac{1}{k(n)}\left| \sum_{j=1}^{k(n)} j \int_{U_{j,n}}^{U_{j+1,n}} \frac{du}{u \log (1/u)} \right|\leq \frac{S^{\ast}_{k(n)}}{k} O_{\mathbb{P}}((\log n)^{-1}).
$$

\Bin where $S^{\ast}_{k(n)}=E_{j,n}+\cdots+E_{k,n}$. We finally get

$$
\sqrt{k(n)} \left(H_n-\gamma\right)=\gamma \frac{S^{\ast}_{k(n)}-k}{\sqrt{k(n)}} + O_{\mathbb{P}}\left(\frac{1}{\log n}, \frac{k^{3/2}}{(\log n)^2}\right)
$$

\Bin We conclude that, whenever (K1) holds, we have

$$
\sqrt{k(n)} \left(H_n-\gamma\right)=\gamma \frac{S^{\ast}_{k(n)}-k}{\sqrt{n}} + o_{\mathbb{P}}(1). \ \square
$$

\Bin 
\Ni \textbf{Proof of the Corollary \ref{theoHillCoro}}. The proof of Part (b) is the conclusion of the proof of Theorem \ref{theoHill} up to the formula \eqref{step_01}. If \label{Freduced} holds, further steps are dismissed. And we need only  $k(n)/n\rightarrow 0$ to conclude. Let us set

$$
 Z^{\ast}_n=\frac{1}{\sqrt{k(n)}} \sum_{j=1}^{k(n)} j\int_{U_{j,n}}^{U_{j+1,n}} \frac{du}{u \log (1/u)}, \ n\geq 1.
$$

\Ni For the first  part, we already knew that $Z_n^{\ast}=O_{\mathbb{P}}(1/\log n)$. We denoted by $W(1)$ a standard Gaussian random variable. By the classical K\'omlos-M\`ajor-Tusn\`ady (KMT) approximation, we have

$$
\left|\frac{S^{\ast}_k(n)-k(n)}{\sqrt{k(n)}}-\gamma W(1)\right|=O_{\mathbb{P}}\left(\frac{\log k(n)}{\sqrt{k(n)}}\right).
$$

\Bin Straightforward expansions using the different rates of convergence lead to

$$
\frac{\sqrt{k(n)}(H_n-\gamma)-\gamma W(1)}{\gamma Z_n^{\ast}} \rightarrow_{\mathbb{P}} 1,
$$

\Bin whenever  $k(n)/n \rightarrow 0$. Now we apply Proposition  in \cite{ips-wciaa-ang}, page 22. Since the function $\log (1/u)$ is slowly varying and that 
$U_{1,n}/U_{k+1,n}$ and $U_{k+1,n}/U_{1,n}$ are both asymptotically bounded in probability, we have

$$
t_n=\sup_{1\leq j \leq k(n)} \sup_{s \in [U_{j,n},U_{j+1,n}]} \left|\frac{\log (1/s)}{\log n}-1\right| \rightarrow_{\mathbb{P}} 0.
$$

\Bin It comes that

$$
Z_n^{\ast}=\frac{\sqrt{k(n)}}{\log n} (k_{(n)}^{-1} S^{\ast}_{k(n)})(1 + O(t_n))=\frac{\sqrt{k(n)}}{\log n} (1+o(1)),
$$

\Bin which gives the desired result. $\blacksquare$\\

\Bin We have the following convergence of the Double-indexed functional Hill statistics.

\begin{theorem} \label{theoDH} We have the following two results.\\

\Ni (a) If the following conditions hold, as $n\rightarrow +\infty$

$$
s_n(f,1)/(s_n(f,s) \log n)\rightarrow 0 \ and \ B_n(f,s) \rightarrow 0,
$$

\Bin then

$$
\frac{T_n(f,s)-\gamma^s a_n(f,s)}{s_n(f,s)} \rightsquigarrow \mathcal{N}\left(0,\gamma^{2s}\right).
$$

\Bin (b) Furthermore, if $a_n(f,s)/a_n(f,s)\rightarrow +\infty$, then

$$
\frac{s_n(f,s)}{a_n(f,s)} \biggr(\left(\frac{T_n(f,s)}{a_n(f,s)}\right)^{1/s}-\gamma\bigg) \rightsquigarrow \mathcal{N}(0,s^{-2}\gamma^2). \blacksquare.
$$
\end{theorem}

\Bin \textbf{Proof}. Let us exploit the proof of Theorem \ref{theoHill}. we have for $j \in \{1,\cdots,k(n)\}$, $s\geq 1$,
\begin{eqnarray*}
A_{i,n}&=&f(j) \left(X_{n-j+1,n}-X_{n-j,n}\right)^s\\
&=& f(j) \biggr(\gamma j^{-1}E_{j,n} - \gamma \int_{U_{j,n}}^{U_{j+1,n}} \frac{du}{u \log (1/u)}  + O_{\mathbb{P}}\left(F_{k(n)} \left(log n\right)^{-2}\right)\biggr)^{s}\\
&=:& f(j) \left(\gamma j^{-1}E_{j,n} - R_{j,n}+C_{j,n}\right)^{s}, 
\end{eqnarray*}

\Bin with

$$
C_{j,n}= O_{\mathbb{P}}\left(\left(log n\right)^{-2}\right) (uniformly \ in \ j),  \left|\gamma \int_{U_{j,n}}^{U_{j+1,n}} \frac{du}{u \log (1/u)} \right| \leq 
\frac{\gamma j^{-1} E_{j,n}b(n)}{|\log n - Z_n|}.
$$

\Bin We get, by the mean value theorem, $j \in \{1,\cdots,k(n)\}$, $s\geq 1$,

\begin{eqnarray*}
&&A_{i,n}-\gamma^s f(j)j^{-s} E_{j,n}^{s}\\
&&\leq sf(j) \left|R_{j,n}+C_{j,n}\right| \left(\gamma j^{-1} E_{j,n} +\left|R_{j,n}\right|+\left|C_{j,n}\right|\right)^{s-1}\\
&& \leq \biggr(\frac{s \gamma f(j) j^{-1}E_{j,n}}{|\log n -Z_n|}\biggr)\biggr(\gamma j^{-1} E_{j,n} +\left|R_{j,n}\right|+\left|C_{j,n}\right|\biggr)^{s-1}.
\end{eqnarray*}

\Bin In the lines below, we will bound the term with the power $s-1$. But if $s=1$, there will is nothing to bound. So formulas regarding that term are dismissed for $s=1$ and are used only for $s>2$. For $s\geq 1$, we will use the $C_{s-1}$ inequality ( for $s\leq 2$, with $|a+b|^{s-1} \leq 2^{s-2} |a|^{s-1} + |b|^{s-1}$
$C_{s-1}=2^{s-2}$). For $0<r<1$, it can be easily checked that, for $u>0$ fixed,  the function $g(v)=(u+v)^r-u^r-v^r$ of $v\geq 0$ takes the value $g(v)=0$ and has a non-positive derivative function, so that $g(v)\leq g(0)=0$ for any $v\geq 0$, which is equivalent to $(u+v)^r\leq u^r+v^r$. We finally have that $|a+b|^{s-1} \leq D_{s} |a|^{s-1} + |b|^{s-1}$ with $D_s=1$ for $1<s<2$ and $D_s=C_{s-1}$ for $s\geq 2$. Applying that inequality leads, $j \in \{1,\cdots,k(n)\}$, $s\geq 1$, to

\begin{eqnarray*} 
&&A_{i,n}-\gamma^s f(j)j^{-s} E_{j,n}^{s}\ \ \ \ \ \ \ \ \ \ \ \ \ \  \ \ \ \ \ \ \ \ \ \ \ \ \ \ \ \ \ \ \ \ \ \ \ \ \ \ \  \ \ \ \ \ \ \ \ \ \ \ \ \  \ \ \ \ \ \ \ \ \ \ \  \ \ \ \ \ \ \ \ \ \ \ \ \ (A)\\
&&\leq \biggr(\frac{s \gamma f(j) j^{-1}E_{j,n}}{|\log n -Z_n|}\biggr)\biggr(D_s\gamma^{s-1} j^{s-1}E_{j,n}^{s-1}+\frac{D_{s}^2\gamma^{s-1}j^{s-1}E_{j,n}^{s-1}}{(|\log n - X_n|^{s-1})} +  O_{\mathbb{P}}\left(\frac{D_s^2}{(\log n)^{2(s-1)}}\right)\biggr).
\end{eqnarray*}

\Bin Let us denote

$$
S_n(f,s)=\sum_{j=1}^{k(n)} f(j) j^{-s} E_{j,n}^s
$$

\Bin and

$$
T_n(f,s)=\sum_{j=1}^{k(n)} f(j) \left(X_{n-j+1,n}-X_{n-j,n}\right)^s,
$$

\Bin By combining the results above, we arrive at \label{B}

\begin{eqnarray*} 
&&\biggr|T_n(f,s)-\gamma^s S_n(f,s)\biggr| \ \ \ \ \ \ \ \ \ \ \ \ \ \  \ \ \ \ \ \ \ \ \ \ \ \ \ \ \ \ \ \ \ \ \ \ \ \ \ \ \  \ \ \ \ \ \ \ \ \ \ \ \ \  \ \ \ \ \ \ \ \ \ \ \  \ \ \ \ \ \ \ \ \ \ \ \ \ (B)\\
&&\leq \biggr(\frac{s \gamma S_n(f,1)}{|\log n -Z_n|}\biggr)\biggr(D_s\gamma^{s-1} S_n(Id,s-1)+\frac{D_{s}^2\gamma^{s-1}S_n(Id,s-1)}{(|\log n - Z_n|^{s-1})} +  
O_{\mathbb{P}}\left(\frac{D_s^2}{(\log n)^{2(s-1)}}\right)\biggr).
 \end{eqnarray*}

\Bin Let us study  $S_n(f,s)$. As a sequence of partial sums of real-value independent random variables indexed by $j \in \{1,\cdots,k(n)\}$ with first and second moments

$$
\Gamma(s+1)f(j)j^{-s} \ and \  (\Gamma(2s+1)-\Gamma(s+1)^2) f(j)^2j^{-2s},
$$  

\Bin the asymptotic normality is given by the the theorem of Levy-Feller-Linderberg (See Theorem 20 in \cite{ips-mfpt-ang} we apply to the centered \textit{rrv}'s $\xi_j=f(j) j^{-s}(E_{j,n}^s-\Gamma(s+1))$, after remarking that

$$
\left\{\frac{\mathbb{V}ar(\xi_j) }{\sum_{j=1}^{k(n)} \mathbb{V}ar(\xi_j)}, \ 1\leq j\leq k(n)\right\}=C(s) B_n(f,s).
$$

\Bin So,  as $n\rightarrow +\infty$,

$$
\left(\frac{1}{s_n(f,s)} \sum_{j=1}^{k(n)} \left(f(j) j^{-s} (E_{j,n}^s-\Gamma(s+1))\right) \rightsquigarrow \mathcal{N}(0,1)\right) \ and \ B_n(f,s) \rightarrow 0
$$

\Bin and the Lynderberg condition holds, that is, for any $\varepsilon>0$,

$$
g(n,\varepsilon)=\frac{1}{s_n(f,s)} \sum_{j=1}^{k(n)} \int_{(|\xi_j|>\varepsilon s_n(f,s))} \xi_j^2 \ d\mathbb{P} \rightarrow 0.
$$

\Bin But, for $K^2(s)=\Gamma(4s+1)-4\Gamma(3s+1)\Gamma(s+1)+\Gamma(2s+1)\Gamma(s+1)2-3\Gamma(3s+1)^4$,

$$
\mathbb{E} \xi^4 = K(s) f(s)^4 j^{-4s}
$$

\Bin and, by the Cauchy-Schwarz inequality

\begin{eqnarray*}
&&\int_{(|\xi_j|>\varepsilon s_n(f,s))} \xi_j^2 \ d\mathbb{P}\\
&&\leq \left(\int \xi_j^4 \ d\mathbb{P}\right)^{1/2} \left(\int 1_{(|\xi_j|>\varepsilon s_n(f,s))} \ d\mathbb{P}\right)^{1/2}\\
&& =K f(j)^2 j^{2s} \left(\int 1_{(|\xi_j|>\varepsilon s_n(f,s))} \ d\mathbb{P}\right)^{1/2}\\
&& =K f(j)^2 j^{2s} \mathbb{P}\left(|\xi_j|>\varepsilon s_n(f,s)) \right)^{1/2}\\
&& \leq K f(j)^2 j^{2s} \left(\frac{K(s)^2 f(j)^4 j^{-4s}}{\varepsilon^4 s_n^4(f,s)}\right)^{1/2}\\
&& =K(s)^2 \left(f(j)^2 j^{2s}\right)^2  (s_n^{-2}(f,s)^2\\
&& =\frac{C(s)}{K(s)}B_n(f,s) \frac{\mathbb{V}ar(\xi_j)}{s^{2}\left(f,s\right)}\\
\end{eqnarray*}

\noindent So

$$
g(n,\varepsilon)= \left(\frac{K(s)}{C(s)}\right)^{2}B_n(f,s) \rightarrow 0.
$$

\Bin Our hypothesis $B_n(f,s)\rightarrow 0$ makes the Lynderberg hold and the central limit theorem holds for $S_n(f,s)$, that is

$$
\frac{S_n(f,s)-\gamma^{s}a_n(f,s)}{s_n(f,s)} \rightsquigarrow \mathcal{N}(0,1)
$$

\Bin Now, let us return to the approximation (B) at page \pageref{B}. We have that for $s=1$, the expression denoted as $C_n$ between the pair of big parentheses should be equal to one as explained before. If $s>1$, we have $\sigma^2(s)=\sum_{j\geq 1} j^{-2(s-1)}<+\infty$, we apply a theorem of Kolmogorov (see \cite{ips-mfpt-ang}, Proposition 25, page 233), $S_n(Id,s-1)$ weakly converges to the random variable $W(s)$ with variance $\sigma^2(s)$. Hence $C_n=O_{\mathbb{P}}(1)$. We arrive at

\begin{equation} 
\biggr|\frac{T_n(f,s)-a_n(f,s)}{s_n(f,s)}-\frac{\gamma^{s}(S_n(f,s)-a_n(f,s))}{s_n(f,s)}\biggr|\leq  O_{\mathbb{P}}\bigg(\frac{S_n(f,1)}{s_n(f,s)\log n}\biggr).
 \end{equation}

\Bin The later bound goes to zero in probability if and only if $s_n(f,1)/(s_n(f,s) \log n)\rightarrow 0$. Now, we have

$$
\frac{s_n(f,s)}{a_n(f,s)} \biggr(\frac{T_n(f,s)}{a_n(f,s)}-\gamma^s\biggr)=Z_n+ o_{\mathbb{P}}(1).
$$

\Bin If $s_n(f,s)/a_n(f,s)\rightarrow +\infty$, we can use the $\delta$-method applied to $g(t)=t^{1/s}$ to get

$$
\frac{a_n(f,s)}{s_n(f,s)} \biggr(\left(\frac{T_n(f,s)}{a_n(f,s)}\right)^{1/s}-\gamma\biggr) \rightsquigarrow \mathcal{N}(0,s^{-2}\gamma^2). \blacksquare.
$$

\Bin \textbf{Remark}. In \cite{ngomlo2016}, we gave a direct proof of the asymptotic normality of $S_n(f,s)$ by using the two hypotheses $B_n(f,s)\rightarrow 0$ and
$s_n(f,s) \rightarrow +\infty$. Here, it seems that we only used the first one. But that one could not hold if $S_n(f,s)$ contains a sub-sequence converging to a finite and positive number. That remark should be recalled in interpreting the results in \cite{ngomlo2016}.\\

\section{Upper records values}

\noindent The main result is :

\begin{theorem} If, for each $n\geq 1$, $X^{(n)}$ stands for $n$-th record value, we have as $n\rightarrow +\infty$, 
$$
\frac{X^{(n)}-\gamma n}{\gamma \sqrt{n}} \rightsquigarrow \mathcal{N}(0,1).
$$
\end{theorem}

\Bin \textbf{Remark}. We refer the reader to  \cite{lo2019ALR} for a simple introduction to records theory.\\

\noindent \textbf{Proof}. We already noticed that $Z=exp(X)$ is the extremal domain of attraction of $G_{\gamma}(x)=\exp(-(1+\gamma x))$, for $\gamma x>-1$. From Part (b) of Theorem 1 in \cite{lo2019ALR}, the $n$-th record $Z^{(n)}=exp(X^{(n)})$ have the representation

\begin{equation} 
\left(\frac{\exp(X^{(n)})}{H^{-1}(1-e^{-n}}\right)^{1/\sqrt{n}}=\exp(\gamma S_{n}^{\ast}) +o_{\mathbb{P}}(1)
\end{equation}

\noindent where $S_{n}^{\ast}$ has the same law as $\gamma^{-1}(T_n-n)/\sqrt{n}$ with $T_n$ denoting a $\gamma$ law with parameters $n$ and $1$. Since 
$H^{-1}(1-u)=\exp(F^{-1}(1-u))$, we have

\begin{equation} 
\frac{X^{(n)}-F^{-1}\left(1-e^{-n}\right)}{\gamma \sqrt{n}}=S_{n}^{\ast} + o_{\mathbb{P}}(1)
\end{equation}

\noindent By the central limit theorem, it comes that

\begin{equation} 
\frac{X^{(n)}-F^{-1}\left(1-e^{-n}\right)}{\gamma \sqrt{n}}=\mathcal{N}(0,1) + o_{\mathbb{P}}(1).
\end{equation}

\Bin By using Formula \eqref{quantile}, we get

\begin{eqnarray} 
&&\frac{X^{(n)}-\gamma n}{\gamma \sqrt{n}}=S^{\ast}_n + o_{\mathbb{P}}(1) \label{records01}\\
&&\frac{X^{(n)}-\gamma n}{\gamma \sqrt{n}}=\mathcal{N}(0,1) + o_{\mathbb{P}}(1). \label{records02}
\end{eqnarray}

\Bin The proof is over. $\square$

\section{The moment problem}

\Ni Typically, the moment problem on $\mathbb{R}$(see \cite{shohat}) is the following. Given a sequences real numbers $(m_n)_{n\geq 1}$, can we find a distribution 
(not necessarily a \textit{cdf}) $F$ on $\mathbb{R}$ as the unique solution of the moments equations.

$$
\forall n\geq 1, \ m_n=\int x^n \ dF(x). 
$$

\Bin This is a nice but difficult mathematical question treated in \cite{shohat}. But in the context of probability theory on $\mathbb{R}$, we may have a fixed \textit{cdf} $F$ of random variable $X$ having moments

$$
\forall n\geq 1, \ \mathbb{E} X^n = m_n \ finite.
$$ 

\noindent The moment problem becomes : Is the sequence of moments $(m_n)_{n\geq 1}$ characterize the probability law of $X$. In that regard, we have

\begin{theorem} \label{momProb} The moments of the pseudo-Lindely probability law are the following

$$
\forall n\geq 1, \ m_n=\frac{n! (\beta+n)}{\theta^n \beta}.
$$

\Bin Any real-valued random variable have the moments $(m_n)_{n\geq 1}$ follows the pseudo-Limdley law.  
\end{theorem}

\textbf{Proof}. At the place of a simple proof, we proceed to slight round-up of the moment problem and explain how to find a simple criteria based on Analysis. A possible tool is the characteristic function which characterize its associated probability law. We have the following expansion of any characteristic function of $X$ (see \cite{loeve} or \cite{ips-mfpt-ang}, Lemma 5, page 255), we have 

\begin{equation}
\mathbb{E}e^{iuX}=1+\sum_{k=1}^{n}\frac{(iu)^{k}mk}{k!}+\theta 2^{1-\delta }\mu
^{n+\delta }\frac{\left\vert u\right\vert ^{n+\delta }}{(n+1)!}.  \label{expan-ch}
\end{equation}

\Bin By usual analysis tools, the series in Formula \eqref{expan-ch} converges in the $]-R, R[$ where $R$ is found according the Cauchy rule

$$
\limsup_{n\rightarrow +\infty} (m_n)^{1/n} = R > 0.
$$

\Bin The conclusion is that two random variables have the same moments of all orders have characteristic functions coinciding on $]-R, R[$. Finally, (see  \cite{loeve}, page 225, Part B.; see also \cite{billinsgleypm}) two characteristic function coinciding on an interval $]-R, R[$ coincide everywhere and thus, are associated to the same probability law.\\

\Bin Let us apply to the pseudo-Lindley law. In \cite{Zeghoudi3}, the moments are given by

$$
\forall n\geq 1, \ m_n=\frac{n! (\beta+n)}{\theta^n \beta}.
$$

\Bin Straightforward computation based on the Stirling formula leads to $R=1/\theta$. This is enough to prove the claim of the theorem.


\newpage
\Ni \textbf{Appendix} \label{appendix-erZ}. Let $R=\beta/\theta$. In the computations below, $u \in (0,1)$ and $x\geq 0$ are linked by  $u=1-F(x)$. So $u\rightarrow 0$ if and only if $x\rightarrow +\infty$. Also, below, functions of $x$ are functions of $u$ actually. We denote $A(u)=\log(1+R/x)$. We have $A(u)\rightarrow 0$ as $u\rightarrow 0$. By writing

$$
\log(\beta +\theta x)=\log(\beta +\theta x)-\log \theta x+\log \theta x= \log \theta x +A(u),
$$

\Bin we see that $u=1-F(x)$ gives

\begin{equation} 
\theta x = \log(1/u) +\log R + \log x + A(u). \label{erZ_01}
\end{equation}

\Bin So, we have

\begin{equation} \label{erZ_02}
F^{-1}(1-u)=\theta^{-1} \log(1/u) (1+o(1)).
\end{equation}

\Bin and

\begin{equation} \label{erZ_03}
\log x = \log\log(1/u) (1+o(1)).
\end{equation}

\Bin Now, we wish to develop that asymptotic equivalence with rates of convergence. Let $B(u)=\log R + \log x+ A(u)$. From Formula \ref{erZ_01}, we have

\begin{equation} \label{erZ_04}
\frac{x}{\theta^{-1}\log(1/u)}-1=\frac{B(u)}{\log(1/u)}.
\end{equation}

\Bin By Formula \eqref{erZ_04}, we notice that 

\begin{equation} \label{erZ_05}
B(u)= \log R + \log x + (R/x) - (R/x)^2/2 + O(\log(1/u)^{-3})=O(\log x)= (\log\log u)(1+o(1)),
\end{equation}

\Bin and hence, for $D(u)=\log R + A(u)$, 

\begin{equation} \label{erZ_06}
\frac{\log(1/u)}{\log x}\left(\frac{x}{\theta^{-1}\log(1/u)}-1\right)= 1 + \frac{D(u)}{\log x}.
\end{equation}

\noindent Also

$$
\frac{D(u)}{\log x}=\frac{\log R + (R/x) - (R/x)^2/2 + O(x^{-3}}{\log x}
$$

\Bin Next, we  have

\begin{eqnarray} \label{erZ_07}
&&\frac{\log x}{-\log R}\left(\frac{\log(1/u)}{\log x}\left(\frac{x}{\theta^{-1}\log(1/u)}-1\right)-1\right)\\
&&=1+\frac{R}{x \log R} -\frac{R^2}{2x^2\log R}+O(x^{-3}) \notag
\end{eqnarray}

\Bin and finally

\begin{eqnarray} \label{erZ_08}
&&\frac{x \log R}{R}\left(\frac{\log x}{\log R}\left(\frac{\log(1/u)}{\log x}\left(\frac{x}{\theta^{-1}\log(1/u)}-1\right)-1\right)-1\right)\\
&&=1-\frac{R}{2x}+O(x^{-2}).\notag
\end{eqnarray}

\Bin Now we want to do the same for $\log x$. From Formula \label{erZ_04}, we have

\begin{equation} \label{erZ_09}
\log (\theta x)= \log\log(1/u) + \log(1+ B(u)/\log(1/u))
\end{equation}

\Bin from which we get

\begin{equation} \label{erZ_09}
\log x - \log\log(1/u)=-\log \theta  + (B(u)/\log(1/u))+ O\left((B(u)/\log(1/u)^2\right).
\end{equation}

\Bin From Formula \eqref{erZ_06}, we have 

\begin{eqnarray*} 
&&\frac{\log(1/u)}{\log x}\left(\frac{x}{\theta^{-1}\log(1/u)}-1\right)-\frac{\log(1/u)}{\log\log 1/u}\left(\frac{x}{\theta^{-1}\log(1/u)}-1\right)\\
&&=\left(\frac{x}{\theta^{-1}\log(1/u)}-1\right)\frac{-(\log(1/u))(\log x -\log\log 1/u)}{(\log x)(\log\log 1/u)}\\
&&=(1+D(u)/\log x)\left(\frac{1}{(\log x)(\log\log 1/u)}\left( -\log \theta  + (B(u)/\log(1/u))+ O(B(u)/\log(1/u)^2)   \right)\right)\\
&&=O((\log\log 1/u)^{2})
\end{eqnarray*}

\noindent Formula \eqref{erZ_06} becomes

\begin{equation} \label{erZ_10}
\frac{\log(1/u)}{\log\log 1/u}\left(\frac{x}{\theta^{-1}\log(1/u)}-1\right)= 1 + \frac{D(u)}{\log x}+ O((\log\log 1/u)^{2}).
\end{equation}

\Bin That formula will be used with Formula \ref{erZ_09} and

\begin{eqnarray} \label{erZ_11}
\frac{B(u)}{\log 1/u}&=&\frac{\log R}{\log 1/u} + \frac{\log\log 1/u}{\log 1/u}(1+o(1))\\
&+& \frac{(R/x)-(R/x)^2/2}{\log 1/u} + O((\log 1/u)^{-4}). \notag
\end{eqnarray}
 
\Bin From \ref{erZ_01}, and from the following formula we can check by using differentiation methods to establish monotonicity

$$
x - x^2/2 \leq \log(1+x) \leq x
$$ 

\Bin we have

\begin{equation} \label{erZ_12a}
(R/x)-R^2/(2x^2) +\log R +\log x \leq \theta x - \log(1/u) \leq  (R/x) +\log R +\log x.
\end{equation}

\Bin But we also have

$$
x= \log(1/u) \left(1 + \frac{\log \beta^{-1}+ \log x + A(u)}{\log(1/u)}\right) 
$$

\Bin which implies

$$
\log x= \log\log(1/u) + \log \left(1 + \frac{\log \beta^{-1}+ \log x + A(u)}{\log(1/u)}\right) 
$$

\Bin By putting

$$
H(u)=\frac{\log \beta^{-1}+ \log x + A(u)}{\log(1/u)},
$$

\Bin we finally get

\begin{equation} \label{erZ_12b}
 H(u)-H(u)^2/2  \leq \log x  -\log\log(1/u)\leq H(u). 
\end{equation}

\Bin By combining Formulas \eqref{erZ_12a} and \eqref{erZ_12b}, we get

\begin{equation} \label{erZ_12c}
 \left|\theta x - \log(1/u) - \log(1/u) \right|\leq \frac{1}{2}\left(\frac{R^2}{x^2}+H(u)^2\right).
\end{equation}

\Bin Since $(R/x^2)$ and $H(u)^2$ are both $O(\log 1/u)^{-2})$, we have

\begin{equation} \label{erZ_12d}
 F^{-1}(1-u) = \theta^{-1}(\log(1/u) -  \log\log(1/u)) +  O(\log 1/u)^{-2}).
\end{equation}

\Bin But since the derivative $\log\log(1/u)$ is $(-u \log(1/u))^{-1}$, we have for $d=-\log\log 2$,

$$
\forall u \in ]0,1[, \ \log\log(1/u)- =\int_{u}^{1/2} \frac{1}{u \log(1/u} \ du,
$$

\Bin and finally

\begin{equation} \label{erZ_12e}
 F^{-1}(1-u) = d + \theta^{-1}(\log(1/u) - \int_{u}^{1/2} \frac{1}{u \log(1/u)} \ du +  O\left(\left(\log 1/u\right)^{-2}\right).
\end{equation}



\end{document}